\documentclass[12pt,a4paper]{article}
\newenvironment{proof}{\bigskip\noindent\emph{Proof.}}{\bigbreak}

\usepackage{amscd,amsmath,amssymb,amsfonts}

\newcommand{\IA}{{\mathbb A}}
\newcommand{\IP}{{\mathbb P}}

\newcommand{\Spec}{\operatorname{Spec}}
\newcommand{\Proj}{\operatorname{Proj}}
\newcommand{\Sym}{\operatorname{Sym}}

\newcommand{\fm}{\mathfrak{m}}
\newcommand{\fn}{\mathfrak{n}}

\newcommand{\noi}{\noindent}
\newcommand{\bs}{\bigskip}
\newcommand{\ms}{\medskip}

\def\qss{quasi-semistabe }
\def\cQS{{\mathcal Q}{\mathcal S}}
\def\cS{{\mathcal S}}
\def\cC{{\mathcal C}}
\def\cO{{\mathcal O}}
\def\Xsred{X_{s,red}}

\def\ch{\mathrm{ch}}

\def\Coker{\mathrm{Coker}}

\def\recmap#1{\rho_{#1}}
\def\recmapn#1{\rho_{#1,n}}

\def\recmapln#1{\rho_{#1,\ell^\nu}}

\def\scs{\; : \;}
\def\piab#1{\pi_1^{ab}(#1)}
\def\sumd#1#2{\underset{x\in {#1}_{#2}}\bigoplus}

\def\rmapo#1{\overset{#1}{\longrightarrow}}

\def\BKl#1#2{BK_{#1}(#2,\ell)}
\def\BKtwo#1#2{BK_{#1}(#2,2)}

\def\tX{\widetilde{X}}

\def\bP{{\mathbb P}}

\def\VS{V_{\Sigma}}

\def\Claim#1{\bigskip\noindent\textbf{Claim #1}\; }

\def\CS{C_{\Sigma}}

\def\isom{\overset{\cong}{\longrightarrow}}

\def\nz{\Bbb Z/n\Bbb Z}

\begin{document}

\begin{center}
{\bf\large Bertini theorems and Lefschetz pencils over discrete
valuation rings, with applications to higher class field theory}
\end{center}

\begin{center}{Uwe Jannsen and Shuji Saito}\end{center}

\bigskip
Good hyperplane sections, whose existence is assured by Bertini's
theorem, and good families of hyperplane sections, so-called
Lefschetz pencils, are well-known constructions and powerful
tools in classical geometry, i.e., for varieties over a field. But
for arithmetic questions one is naturally led to the consideration
of models over Dedekind rings and, for local questions, to schemes
over discrete valuation rings. It is the aim of this note to
provide extensions of the mentioned constructions to the latter
situation. We point out some new phenomena, and give an
application to the class field theory of varieties over local
fields with good reduction or, more generally, ordinary quadratic reduction.
For this we also discuss the resolution of the latter to semi-stable models.
For more arithmetic applications see \cite{JS}.

Let $A$ be a discrete valuation ring with fraction field $K$,
maximal ideal $\mathfrak{m}$ and residue field $F = A/\fm$. Let
$\eta = \Spec(K)$ and $s = \Spec(F)$ be the generic and closed
point of $\Spec(A)$, respectively. For any scheme $X$ over $A$ we
let $X_\eta = X\times_AK$ and $X_s = X\times_AF$ be its generic
and special fibre, respectively.

\bs \setcounter{section}{0}
\begin{center}{\bf\arabic{section}. Good hyperplane sections for good reduction schemes}\end{center}

As a `warm-up', we recall the classical Bertini theorem and extend
it to varieties over $K$ with good reduction. Let
$X\subset\IP^N_L$ be a smooth quasi-projective variety over a
field $L$. Recall that another irreducible smooth subscheme
$Z\subset \IP^N_L$ is said to intersect $X$ transversally, if the
scheme-theoretic intersection $X\cdot Z = X\times_{\IP^N}Z$ (which
is just defined by the ideal generated by the equations of $X$ and
$Z$) is smooth and of pure codimension $\mathrm{codim}_{\IP^N}(Z)$
in $X$. Then the Bertini theorem asserts that for infinite $L$,
there exists an $L$-rational hyperplane $H\subset \IP^N_L$
intersecting $X$ transversally (cf. [Jou, 6.11, 2], and also
Theorem 3 below). In this case, one calls $Y = X\cdot H$ a smooth
(or good) hyperplane section of $X$.

More precisely, the following holds. One has the dual projective
space $(\IP^N_L)^\vee$ parameterizing the hyperplanes in $\IP^N_L$
(a point $a = (a_0:\ldots:a_N)$ corresponds to the hyperplane with
the equation $a_0x_0 + \ldots + a_Nx_N = 0$ for the homogeneous
coordinates $x_i$ of $\IP^N_L$). Then, for an arbitrary field $L$,
there is a dense Zariski open $V_X\subset (\IP^N_L)^\vee$
parameterizing those hyperplanes which intersect $X$
transversally. Moreover, if $L$ is infinite, then the set $V_X(L)$
of $L$-rational points is non-empty, since $\IP^N_L(L)$ is Zariski
dense in $\IP^N_L$. This shows that, for an infinite field $L$,
and finitely many smooth varieties $X_1,\ldots,X_n$ in $\IP^N_L$,
there also exists an $L$-rational hyperplane $H$ intersecting all
$X_i$ transversally, because $V_{X_1}\cap\ldots\cap V_{X_n}$ is
non-empty.

If $L$ is finite, it may happen that $V_X$ does not have any
$L$-rational point. But, by sieve methods, Poonen [Po] showed that
in this case there always exists an $L$-rational point after
replacing the projective embedding by the $d$-fold embedding for
some $d > 0$, i.e., there always exists a smooth $L$-rational {\it
hypersurface} section of $X$.

Now consider a quasi-projective $A$-scheme ($A$ a discrete
valuation ring as above), i.e., a subscheme $X$ of the projective
space $\IP^N_A$ over $A$.

By a hyperplane $H\subseteq\IP^N_A$ over $A$ we mean a closed
subscheme which corresponds to an $A$-rational point of the dual
projective space $(\IP^N_A)^\vee$ (= Grassmannian of linear
subspaces of codimension $1$). Since every invertible module over
$A$ is free, $H$ is given by a surjection
$\varphi:A^{N+1}\twoheadrightarrow A^N$; or, equivalently, by an
equation $\sum^N_{i=0}a_ix_i=0$, $a_i\in A$ $(i=0,\dots, N)$, not
all in the maximal ideal $\fm$, for the coordinates $x_i$ on
$\IP^n_A$. The correspondence is given by
$$
\ker\varphi=A\cdot \sum^N_{i=0}a_ie_i\qquad,
$$
where $e_0,\dots, e_N$ is the basis of $A^{N+1}$.

\ms\noi {\bf Theorem 0} ~{\it Let $X\subset \IP^N_A$ be a smooth
quasi-projective $A$-scheme. If $F$ is infinite, then there exists
a hyperplane $H\subset \IP^N_A$ over $A$ such that the
scheme-theoretic intersection $X\cdot H = X\times_{\IP^N_A}H$ is
smooth over $A$ and of pure codimension 1 in $X$. If $F$ is finite
and $A$ is Henselian, then, for every given prime number $\ell$,
such a hyperplane exists after replacing $A$ by a finite \'etale
covering $A'/A$ of $\ell$-power-degree.}

\ms\noi {\bf Proof} ~Let $H\subset \IP^N_A$ be a hyperplane over
$A$. Then $H_\eta$ and $H_s$ are hyperplanes in $\IP^N_K$ and
$\IP^N_F$, respectively. With the notations as above, the
condition on the hyperplane is that (the $K$-rational point
corresponding to) $H_\eta$ lies in the good locus
$V_{X_\eta}\subset (\IP^N_K)^\vee$, and that $H_s$ lies in
$V_{X_s}$. Since $H$ is completely determined by $H_\eta$, this
means that $H_\eta \in V_{X_\eta}(K)\cap sp^{-1}(V_{X_s}(F))$,
where $sp: (\IP^N)^\vee(K) \rightarrow (\IP^N(F)$ is the
specialization map, which sends $H_\eta$ to $H_s$.

It remains to see when this intersection is non-empty. But for any
proper scheme $P$ over $A$ and any open subschemes $V_1\subset
P_\eta$ and $V_2\subset P_s$, with closed complements $Z_1 =
P_\eta \setminus V_1$ and $Z_2 = P_s\setminus V_2$, respectively,
one has $Z_1(K)\subset sp^{-1}(sp(Z_1)(F))$ where $sp(Z_1) =
\overline{Z}_1\cap P_s$ for the Zariski closure $\overline{Z}_1$
of $Z_1$ in $\IP^N_A$. Therefore $V_1(K)\cap sp^{-1}(V_2(F))$
contains $sp^{-1}((V_2\setminus sp(Z_1))(F))$. The latter set has
$K$-rational points, if $sp: P(K) \rightarrow P(F)$ is surjective
and $V_2\setminus sp(Z_1)$ has $F$-rational points. The latter set
is open and dense in $P_s$, if $P/S$ has irreducible fibres, and
$V_1$ and $V_2$ are dense in their fibres.

Applying this to $P=(\IP^N_A)^\vee$, $V_1=V_{X_\eta}$ and
$V_2=V_{X_s}$, where all conditions are fulfilled, we see it
suffices that the non-empty open subset $W=V_2\setminus sp(Z_1)$
has $F$-rational points. As explained above, this is the case if
$F$ is infinite. Hence, if $F$ is finite, it is the case over the
maximal pro-$\ell$-extension of $F$, hence over some extension
$F'/F$ of $\ell$-power degree. If $A$ is Henselian, and $A'/A$ is
the unramified extension corresponding to $F'/F$, then the
$F'$-rational point lifts to an $A'$-rational point of $P$. Since
the formation of the sets $V_1$ and $V_2$ is compatible with
\'etale base change in the base, this means there is a good
hyperplane section for $X$ over $A'$.

\bs\noi {\bf Remarks 0} ~(i) In contrast with the classical
situation, the good hyperplanes over $A$ are not parametrized by a
Zariski open in $\IP^N_K$, but by a subset of the type $V_1(K)\cap
sp^{-1}((V_2)(F))$ for Zariski opens $V_1\subset (\IP^N_K)^\vee$
and $V_2\subset (\IP^N_F)^\vee$.

(ii) If, with the notations as in the proof, $H_s$ intersects the
smooth variety $X_s$ transversally, and if $X_\eta\cap H_\eta$ is
non-empty, then $X\cdot H$ is a flat $A$-scheme of finite type,
whose special fibre $(X\cdot H)_s = X_s\cdot H_s$ is smooth. Since
the smooth locus of $X\cdot H$ is open, $X\cdot H$ must be smooth,
if $X$ and hence $X\cdot H$ is proper. This shows that, for smooth
and proper $X$, one has $sp^{-1}(V_{X_s}) \subset V_{X_\eta}$, and
the locus in $(\IP^N_K)^\vee(K)$ of good hyperplanes for $X/A$ is
just $sp^{-1}(V_{X_s}(F))$. Moreover, by applying the mentioned
result of Poonen, this has a $K$-rational point after passing to
some multiple embedding.

\bs Recall that a smooth proper variety $V$ over $K$ is said to
have good reduction (over $A$) if there is a smooth proper
$A$-scheme $X$ with generic fiber $X_\eta = X\times_AK \cong V$.

\bs\noi {\bf Corollary 0} ~{\it If $F$ is finite and $V/K$ is a
smooth projective variety with good reduction, there exists a
smooth hypersurface section which again has good reduction.}

\bs\setcounter{section}{1}
\begin{center}{\bf\arabic{section}. Good hyperplane sections
for quasi-semi-stable schemes}\end{center}

For the applications, the case of good reduction is too restrictive.
We now introduce the following more general objects:

\ms\noi {\bf Definition 1}
Let $\cC$ be the category of quasi-projective schemes over $A$.
An object $X\in \cC$ is {\it \qss} if the following conditions holds:
\begin{enumerate}
\item[(1)] $X$ is regular flat over $\Spec(A)$,

\item[(2)] for each closed point $x \in X_s$ the completion of
$\cO_{X,x}$ is isomorphic to
$$
B=A[[x_1,\dots, x_r, y_1,\dots, y_n]]/ \langle \pi- u
x_1^{e_1}\dots x_r^{e_r} \rangle \qquad
$$
where $e_1,\dots,e_r\geq 1$ are suitable integers and $u$ is a
unit in the ring of the formal power series $A[[x_1,\dots, x_r,
y_1,\dots, y_n]]$.

\end{enumerate}
Let $\cQS\subset \cC$ be the subcategory of all $X\in Ob(\cC)$
which are \qss and let $\cS\subset \cQS$ the subcategory of all $X\in
Ob(\cQS)$ for which $X_s$ is reduced. (This means that
$e_1=\cdots e_r=1$ in condition (2), so that $X$ has so-called
{\it semi-stable} reduction.)
\bigskip

The aim of this section is to prove:

\ms\noi {\bf Theorem 1} ~{\it
Let $X$ be an object of $\cQS$ (resp. $\cS$).
If $F$ is infinite, then there
exists a hyperplane $H \subset \IP^N_A$ over $A$ such that
$X$ and $H$ intersect transversally in $\IP^N_A$,
$X \cdot H :=X\times_{\IP^N_A}H$ is in $\cQS$ (resp. $\cS$), and
$(X\cdot H) \cup\Xsred$ is a simple normal crossing divisor on $X$.
If $F$ is finite and $A$ is
Henselian, then for every prime $\ell$, there is a finite
unramified extension $A'$ of $A$ of $\ell$-power degree such that
the same conclusion holds after base change with $A'$.}
\medbreak

For the proof we need the following lemma.

\ms\noi {\bf Lemma 1} ~{\it
Let $X$ be an object of $\cQS$ (resp. $\cS$).
Let $H \subset \IP^N_A$ be a
hyperplane over $A$, with special fibre $H_s\subset\IP^N_F$
and generic fibre $H_\eta \subset \IP^N_K$.
Let $Y_1,\dots,Y_M$ be the irreducible components of $\Xsred$, which are
by definition smooth varieties intersecting transversally in $\IP^N_F$.
Assume that
\begin{itemize}
\item[$(i)$]
$H_s$ and $Y_{i_1,\dots, i_p}$ for any $i_1,\dots,i_a$ intersect transversally
in $\IP^N_F$, where
$Y_{i_1,\dots, i_p}:= Y_{i_1}\cap\dots\cap Y_{i_p}$.
\item[$(ii)$]
$H_\eta$ and $X_\eta$ intersect transversally in $\IP^N_K$.
\end{itemize}
Then $X$ and $H$ intersect transversally in $\IP^N_A$ and
$X \cdot H := X\times_{\IP^N_A}H$ is an object of $\cQS$ (resp. $\cS$)
and $(X\cdot H) \cup\Xsred$ is a simple normal crossing divisor on $X$.
If $X$ is proper over $A$, condition (ii) is implied by condition (i).}
\bigskip

\ms\noi {\bf Proof of Lemma 1} ~Noting
$$
(X\cdot H)\times_X Y_{i_1,\dots, i_p} =
(H\times_{\IP^N_A} X)\times_X Y_{i_1,\dots, i_p} =
H\times_{\IP^N_A} Y_{i_1,\dots, i_p} =
H_s\times_{\IP^N_F} Y_{i_1,\dots, i_p},
$$
it suffices to show that $X \cdot H$ is an object of $\cQS$ (resp.
$\cS$). We may assume that the residue field $F$ of $A$ is
algebraically closed. Choose a closed point $x\in X_s$ and assume
that the completion of $\cO_{X,x}$ is isomorphic to
$$
B=A[[x_1,\dots, x_r, y_1,\dots, y_n]]/
\langle \pi- u x_1^{e_1}\dots x_r^{e_r} \rangle \qquad
$$
as in condition (2) above. Let $f\in B$ be the image of the local
equation for $H$ at $x$, and let $\fn\subseteq B$ be the maximal
ideal. Since $B/\langle f\rangle$ is the completion of the local
ring of $X\cdot H$ at $x$ if $x\in X\cdot H$, and since the
irreducible components of $(X \cdot H)_{s,red} = \Xsred \cap H_s$
are the connected components of the smooth varieties $Y_i \cap
H_s$, the lemma follows from the following two claims. In fact,
Claim 2 shows that every $x\in (X\cdot H)_{s,red}$  has an open
neighborhood in $X\cdot H$ which is an object $\cQS$ (resp.
$\cS$). If $X/A$ is proper, these neighborhoods cover $X\cdot H$.
\bigskip

\ms\noi {\bf Claim 1}  Assumption (i) implies that

\ms\noi
\begin{tabular}{lll}
either  & (a) & $f$ is a unit in $B$, \\
or      & (b) & $s\geq 1$, $f$ is in $\fn$, and has non-zero image in
            $\fn/(\fn^2+\langle x_1,\dots, x_r \rangle)$.
\end{tabular}

\ms\noi
{\bf Claim 2} Assume condition (b) holds. Then
$$
B/\langle f\rangle \cong A[[x_1,\dots, x_r,  z_1,\dots, z_{s-1}]]/
\langle \pi- \overline{u} \cdot x_1^{e_1}\dots x_r^{e_r} \rangle ,
$$
where $\overline{u}$ is a unit of $A[[x_1,\dots, x_r,  z_1,\dots, z_{s-1}]]$.
\medbreak

\medbreak

\ms\noi {\bf Proof of claim 2} The elements $x_i$ and $y_j$ $\mod \fn^2$ form
an $F$-basis of $\fn/\fn^2$ $~(1\leq i\leq r,\; 1\leq j\leq n)$. Hence we have
$$
f=\sum^r_{i=1}a_i x_i +  \sum^s_{j=1}a_{j+r} y_j \quad\mod \fn^2
$$
with elements $a_i,\; a_{j+r}\in A$ which are determined modulo
$\langle\pi\rangle$. If (b) holds, then $a_{j+r} \in A^\times$ for some $j$,
and by possibly renumbering and multiplying $f$ by a unit
we may assume $j=s$ and $a_{r+s}=1$. But then
$$
B/\langle f\rangle \cong A[[x_1,\dots, x_r,  y_1,\dots, y_{s-1}]]/
\langle \pi- \overline{u} \cdot x_1^{e_1}\dots x_r^{e_r} \rangle.
$$

\ms\noi {\bf Proof of claim 1} The elements $x_1,\dots, x_r$ are
the images of the local equations for $Y_{i_1},\dots, Y_{i_r}$ for
suitable $1\leq i_1<\dots <i_r\leq M$. Thus the trace of
$Y_{i_1}\cap\dots\cap Y_{i_r}$ in $\hat{\mathcal O}_{X,x} \cong B$
is given by the ideal $\langle x_1,\dots, x_r\rangle$, i.e., by
the quotient
$$
B'=B/\langle x_1,\dots, x_r\rangle \cong F[[y_1,\dots,y_s]]\quad .
$$
This is zero-dimensional if and only if $s=0$, and in this case
$Y_{i_1}\cap \dots\cap Y_{i_r}$ is zero-dimensional as well. Then,
by assumption on $H$, $H$ does not intersect $Y_{i_1}\cap\dots\cap
Y_{i_r}$, and so $f$ is a unit in $B/\langle x_1,\dots,x_r\rangle$
and hence so is in $B$.

If $s\geq 1$, then $H$ intersects $Y_{i_1}\cap\dots\cap Y_{i_r}$
transversally at $x$ if and only if the image of $f$ in $B'$ lies
in $\fn'-(\fn')^2$, where $\fn'$ is the maximal ideal of $B'$. Now
Claim 1 follows from the isomorphism
$$
\fn'/(\fn')^2\cong \fn/(\fn^2 + \langle x_1,\dots, x_r\rangle).
$$
\medbreak

\ms\noi {\bf Proof of Theorem 1} It suffices to find a hyperplane
satisfying the assumption of Lemma 1, i.e., to show that, with the
notations introduced earlier, the set $V_{X_\eta}(K)\cap
sp^{-1}(V_2(F))$ is non-empty, where $V_2$ is the intersection of
the sets $V_{Y_{i_1,\dots, i_p}}$, and hence open and dense in
$(\IP^N_F)^\vee$. This holds under the conditions of Theorem 1, by
the arguments used in the proof of Theorem 0.
$\square$
\medbreak

\ms\noi If $X/A$ is proper, we noted that $sp^{-1}(V_2(F))$ is
contained in $V_{X_\eta}(K)$. Combining this with the mentioned
results of Poonen [Po] we get:

\ms\noi {\bf Corollary 1} ~{\it If $F$ is finite and $V/K$ is a
smooth projective variety with strictly semi-stable reduction,
there is a smooth hypersurface section which again has strictly
semi-stable reduction.}

\pagebreak
\bigskip
\setcounter{section}{2}
\begin{center}{\bf\arabic{section}. Lefschetz pencils for schemes with
ordinary quadratic reduction}\end{center}

Even if one starts with a variety $V$ over $K$ with good
reduction, in general infinitely many fibres in a Lefschetz pencil
(cf. below) for $V$ will not have good reduction, because
infinitely many hyperplanes specialize to the same hyperplane in
the reduction, and usually the induced pencil for the reduction of
$V$ has a bad member. But one can arrange very mild singularities:

\ms\noi {\bf Definition 2} We say that a smooth projective variety $V$
over $K$ has ordinary quadratic reduction, if there is a
projective $A$-scheme $X$ such that $X_\eta \cong V$, and $X_s$ is
smooth over $F$ except for a finite number of singular points
which are ordinary quadratic (cf. [SGA 7 XV, 1.2.1, XII 1.1] and below).

\ms\noi We will show that one can even start with such singularities, and
still get singularities which are not worse - which is useful for
induction on dimension:

\ms\noi {\bf Theorem 2} ~{\it Let $V$ be projective $K$-variety
with ordinary quadratic reduction, and let $X \subset \IP^N_A$ be a model
of $V$ as in Definition 1. let $d\geq 2$ be an integer and suppose
$F$ is infinite. Then, after possibly passing to the $d$-fold
embedding of $X$, there exists a Lefschetz pencil $\{V_t\}_{t\in
D}$, where $D$ is a line in the dual projective space
$(\IP^N_K)^\vee$, satisfying the following conditions:
\begin{itemize}
\item[(1)] The axis of the pencil has good reduction over $K$.
\item[(2)] There exists a finite subset $\Sigma\subset \IP^1_K$ of
closed points such that for any $t\not\in \Sigma$, $V_t$ has
ordinary quadratic reduction over $A_{K(t)}$, the integral closure of $A$
in the residue field $K(t)$ of $t$.
\end{itemize}
Suppose $F$ is finite, $A$ is Henselian, and $\ell$ is a fixed
prime. Then the same result holds after possibly passing to a
finite unramified extension $K'/K$ of $\ell$-power degree.}

\ms\noi The proof will be achieved in four steps, numbered (2.1)
$\sim$ (2.4).

\ms\noi {\bf (2.1)} Let $(\IP^N_A)^\vee$ be the dual projective
space over $A$. Its fibres over $K$ and $F$ coincide with the dual
projective spaces $(\IP^N_K)^\vee$ and $(\IP^N_F)^\vee$,
respectively. Furthermore, let $\mathcal G=Gr(1,(\IP^N_A)^\vee)$
be the Grassmannian of lines in $(\IP^N_A)^\vee$; again its fibres
over $K$ and $F$ are the corresponding Grassmannians for
$(\IP^N_K)^\vee$ and $(\IP^N_F)^\vee$, respectively. According to
[SGA 7, XVII, 2.5], after possibly passing to the $d$-fold
projective embedding, there is a dense open subscheme $W_{X_\eta}
\subset\mathcal G_K$ such that the lines in $W_{X_\eta}$ give
Lefschetz pencils for $X_\eta \subset \IP^N_K$.

\ms\noi {\bf (2.2)} Since $X_F$ is possibly singular, we need a
slight extension of the results in [SGA 7, XVII]. First we extend
the results to smooth, but only quasi-projective varieties.

\ms\noi {\bf Theorem 3} ~{\it Let $L$ be any field, let $U \subset
\IP^N_L$ be a smooth irreducible quasi-projective variety, and let
$d\geq 2$ be an integer. After possibly passing to the d-fold
embedding, there is a non-empty open subscheme $W_U$ in the
Grassmannian $Gr(1,(\IP^N_L)^\vee)$ of lines in the dual
projective space, such that the lines $D$ in $W_U$ satisfy all
properties of Lefschetz pencils with respect to $U$, i.e.:
\begin{itemize}
\item[(1)] The axis of $D$ (i.e., the intersection of any two
different and hence all hyperplanes parametrized by $D$)
intersects $U$ transversally.
\item[(2)] There is a finite
subscheme $\Sigma \subset D$ such that for $t \in
D\setminus\Sigma$ the hyperplane $H_t$ corresponding to $t$
intersects $U$ transversally.
\item[(3)] For $t \in \Sigma$ the
scheme-theoretic intersection $U \cdot H_t = U\times_{\IP^N_L}H_t$
is smooth except for one singular point which is ordinary
quadratic.
\end{itemize}
}

\ms\noi {\bf Proof} Let $X$ be the closure of $U$ in $\IP =
\IP^N_L$, and let $Z=X\setminus U$ (both endowed with the reduced
subscheme structure). For $Q = \IP^N_L \setminus Z$, let $\mathcal
J$ be the ideal sheaf of the closed immersion $U \subset Q$, and
denote by $\mathcal N = \mathcal J/\mathcal J^2$ the conormal
sheaf, regarded as a locally free sheaf on $U$, and by $\mathcal
N^{\,\vee}$ its dual. As in [SGA 7, XVII] consider the closed
immersion of projective bundles on $U$
$$
\IP_U(\mathcal N^{\,\vee})\hookrightarrow \IP_U(\mathcal
O_U(1)\otimes_L\Gamma(\IP,\mathcal O_\IP(1))^\vee) \cong U\times
(\IP^N_L)^\vee
$$
induced by the canonical monomorphism of bundles
$$
\mathcal J/\mathcal J^2
\hookrightarrow\Omega^1_{Q|U}\hookrightarrow\mathcal O_U(-1)^{N+1}
= \mathcal O_U(-1)\otimes_L\Gamma(\IP,\mathcal O_\IP(1)).
$$
(Here we adopt the convention that, for a vector bundle $\mathcal
F$ on $U$, the projective bundle $\IP(\mathcal F) =
\Proj(\Sym(\mathcal F))$ parametrizes line bundle quotients of
$\mathcal F$.) The above immersion identifies $\IP_U(\mathcal
N^{\,\vee})$ with the subvariety of points $(x,H)$ in $U\times
(\IP^N_L)^\vee$ for which $H$ touches $U$ in $x$. Let $U^\vee$ be
the closure of the image of $\IP_U(\mathcal N^{\,\vee})$ in
$(\IP^N_L)^\vee$. It is the dual variety to $U$ and contains all
hyperplanes in $\IP^N_L$ which touch $U$ in some point. One has
$\dim U^\vee \leq \dim\IP_U(\mathcal N^{\,\vee})=N-1$. Hence
$(\IP^N_L)^\vee \smallsetminus U^\vee$ is non-empty, and the set
$M^{''}_U\subseteq \mathcal G_L = Gr(1,(\IP^N_L)^\vee)$ of lines
in $(\IP^N_L)^\vee$ contained in $U^\vee$ is closed and different
from  $\mathcal G_L$.

Moreover, let $(U^\vee)^0$ be the set of hyperplanes which touch
$U$ in exactly one point which is an ordinary quadratic
singularity. Then $(U^\vee)^0$ is open in $U^\vee$ by results of
Elkik and Deligne ([SGA 7, XVII, 3.2], [SGA 7, XV, 1.3.4]). (If
$\mbox{char } L\neq 2$ or if $n=\dim U$ is even, then it is the
locus where $\IP_U(\mathcal N^{\,\vee})\to (\IP^N_L)^\vee$ is
unramified.) It is non-empty after replacing the given embedding
by its $d$-multiple ($d\geq 2$), by the argument in [SGA 7, XVII,
3.7, 4.2]. Since $U^\vee$ is irreducible (by irreducibility of
$\IP_U(\mathcal N^{\,\vee})$), the closed subscheme $F^{'''}=
U^\vee\smallsetminus(U^\vee)^0$ has codimension $\geq 2$ in
$(\IP^N_L)^\vee$ in this case. Then the set
$M^{'''}_U\subset\mathcal G_L$ of lines in $\IP^N_L$ which meet
$F^{'''}$ is closed and different from $\mathcal G_L$.

Finally, the set $W'_U\subseteq Gr(N-2,\IP^N_L)$ of codimension
$2$ linear subspaces in $\IP^N_L$ which intersect $U$
transversally is open [Jou, 6.11, 2)]. It is also non-empty: Since
$U^\vee \neq (\IP^N_L)^\vee$, over the algebraic closure there is
a hyperplane $H_1$ intersecting $U$ transversally, and similarly,
there is a hyperplane $H_2$ intersecting $U\cdot H_1$
transversally. This means that the codimension $2$ linear subspace
$H_1\cdot H_2$ intersect $U$ transversally. Recall the isomorphism
$$
\mathcal G_L=Gr(1,(\IP^N_L)^\vee)\overset{\sim}\longrightarrow
Gr(N-2,\IP^N_L)
$$
sending a pencil to its axis. We denote the preimage of $W'_U$ in
$\mathcal G_L$ by $W'_U$ again.

The conclusion is that there is a non-empty open subscheme $W_U=
W'_U\cap(\mathcal G_L\smallsetminus M^{''}_U)\cap(\mathcal
G_L\smallsetminus M^{'''}_U)\subseteq\mathcal G_L$ such that the
lines in $W_U$ satisfy all properties of Lefschetz pencils with
respect to $U$, and thus Theorem 3 is proved.

\ms\noi {\bf (2.3)} Now we deal with the singular points of the
special fibre $X_F$ of $X$ in Theorem 2.

\ms\noi {\bf Theorem 4} ~{\it Let $L$ be any field, and let $X
\subset \IP^N_L$ be a projective variety which is smooth except
for finitely many singular points $x_1, \ldots , x_r$ which are
ordinary quadratic. After possibly passing to the d-fold embedding
(any $d\geq 2$), there is a non-empty open subscheme $W_X \subset
Gr(1,(\IP^N_L)^\vee)$ such that for the lines $D$ in $W_X$ the
following holds:
\begin{itemize}
\item[(i)] The axis of $D$ does not meet the singular points of
$X$ and intersects the regular locus $X^{reg}$ transversally.
\item[(ii)] There is a finite subscheme $\Sigma\subseteq D$ such
that for $t\in D\smallsetminus \Sigma$ the hyperplane $H_t$ does
not meet the singular points of $X$ and intersects $X^{reg}$
transversally.
\item[(iii)] For $t\in \Sigma$ the scheme-theoretic
intersection $X^{reg}\cdot H_t$ is smooth except for possibly one
singular point which is an ordinary quadratic singularity.
\item[(iv)] If $t \in \Sigma$ and $x_i \in H_t$, then $x_i$ is an
ordinary quadratic singularity of $X \cdot H_t$.
\end{itemize}
}

\ms\noi {\bf Proof} Applying Theorem 3 to
$X^{reg}=X\smallsetminus\{x_1,\dots, x_r\}$ we find a non-empty
open subset $W'=W_{X^{reg}}\subseteq\mathcal G_L$ such that the
lines $D$ in $V$ satisfies the properties (i) to (iii) for
$X^{reg}$ instead of $X$.

It remains to consider the singular points $x_1,\dots, x_r$. For
each $x_i$, the hyperplanes in $\IP^N_L$ which pass through $x_i$
form a hyperplane $\widetilde H_i\subseteq (\IP^N_L)^\vee$. By the
following lemma there is a non-empty open subset
$U_i\subseteq\widetilde H_i$ such that for any hyperplane $H$ in
$U_i$ the intersection $Y \cdot H$ has an ordinary quadratic
singularity at $x_i$. Then $F_i=\widetilde H_i\smallsetminus U_i$
is closed and of codimension $\geq 2$ in $(\IP^N_L)^\vee$, and so
is $F=\cup^r_{i=1}F_i$. The set $W''\subseteq\mathcal G_L$ of
lines in $(\IP^N_L)^\vee$ which do not meet $F$ and are not
contained in any $\widetilde H_i$ is thus open and non-empty, and
the properties (i) to (iv) above hold for the lines in $W_X =
W'\cap W''\subseteq\mathcal G_L$.

\ms\noi {\bf Lemma 2} ~{\it Let $L$ be any field, let $X \subset
\IP^N_L$ be a projective variety of positive dimension, and let
$x$ be an isolated singularity which is an ordinary quadratic
point. If $\widetilde H_x \subset (\IP^N_L)^\vee$ denotes the
locus of hyperplanes passing through $x$, then there is an open
dense subset $U\subset \widetilde H_x$ such that for all
hyperplanes $H$ in $U$ the point $x$ is an ordinary quadratic
singularity of $X \cdot H$.}

\ms\noi {\bf Proof} We may assume that $L$ is algebraically
closed. Let $A=\hat{\mathcal O}_{X,x}$ be the completion of the
local ring at $x$. Then $x$ is called an ordinary quadratic
singularity, if $A$ is isomorphic to the quotient
$$
L\bigl[[x_1,\dots, x_{n+1}]\bigr]/\langle f\rangle ~,
$$
where $f$ starts in degree $2$, and where $f_2$, the homogeneous
part of degree $2$ of $f$, is non-zero, and defines a non-singular
quadric in $\IP^n_L$ (where $n\geq 1$ by assumption). We shall
call $A$ the ring of an ordinary quadratic singularity in this
case.

\ms\noi {\bf Lemma 3} ~{\it Let $\fm\subset A$ be the maximal
ideal, and let $g\in \fm\smallsetminus\{0\}$ be an element. Then
$A':=A/\langle g\rangle$ is the ring of an ordinary quadratic
singularity if the following two conditions hold
\begin{itemize}
\item[(i)] The image $\bar g$ of $g$ in $\fm/\fm^2$ is non-zero.
\item[(ii)] The non-singular projective quadric
$\Proj(\,\Sym(\fm/\fm^2)/\langle Q\rangle)$ and the hypersurface
$\Proj(\,\Sym(\fm/\fm^2)/\langle\bar g\rangle)$ intersect
transversally in $\Proj(\,\Sym(\fm/\fm^2))\cong \IP^n_L$. Here $Q$
corresponds to $f_2$ under the isomorphism
$$
L[[x_1,\dots, x_{n+1}]]\overset\sim{\longrightarrow}
\Sym(\fm/\fm^2) ~.
$$
(More intrinsically, $Q$ is determined up to a scalar factor as
the generator of the $1$-dimensional kernel of the surjection $
\Sym^2(\fm/\fm^2)\twoheadrightarrow \fm^2/\fm^3 $).
\end{itemize}
}

\ms\noi {\bf Proof} Lift $g$ to an element $\bar g\in
B:=L\bigl[[x_1,\dots, x_{n+1}]\bigr]$, and let $\fn$ be the
maximal ideal of $B$. By (i) and a substitution we may assume
$\bar g=x_{n+1}$. Then $B':=B/\langle\bar g\rangle\cong
L[[x_1,\dots, x_n]]$, and
$$
A'=B'/\langle f'\rangle
$$
where $f'$ is the image of $f$ in $B'$. Then $f'$ starts in degree
$2$ as well, and $f'_2$, its degree $2$ part with respect to the
variables $x_1,\dots, x_n$, is just the image of $f_2$. If $f'_2$
is zero, then $\langle Q\rangle\subseteq\langle\bar g\rangle$ in
$\Sym(\fm/\fm^2)$, in contradiction to (ii). Hence $f'_2\neq 0$,
and by (ii) it gives rise to a non-singular quadric in
$$
\Proj(\,\Sym(\fm/\fm^2)/\langle\bar g\rangle) \cong
\Proj(\,\Sym(\fm/(\fm^2 + \langle g\rangle))\cong \IP^{n-1}_L
$$
by (ii) for $n\geq 2$, i.e., $A'$ is the ring of an ordinary
quadric singularity.

\ms We proceed with the proof of Lemma 2. Choose coordinates
$X_0,\dots, X_N$ on $\IP^N_L$ such that $x=(1:0:\dots:0)$. The
hyperplanes in $\IP^N_L$ are given by points $b=(b_0:\dots:b_N)$
in the dual projective space $(\IP^N_L)^\vee$, corresponding to
the hyperplanes
$$
\sum^N_{i=0}b_iX_i=0\qquad .
$$
The hyperplanes through $x$ are given by those $b$ with $b_0=0$
and are parametrized by $(b_1:\dots:b_N)\subseteq
(\IP^{N-1}_L)^\vee$. If $x_i=\frac{X_i}{X_0}$, $i=1,\dots, N$, are
the affine coordinates on the open affine neighborhood $\{
x_0\neq 0\} \cong \IA^N_L$, $x$ corresponds to the zero point, and
the hyperplane associated to $(b_1:\dots:b_N)$ is determined by
the element $\sum\limits ^N_{i=1}b_ix_i\in L[x_1,\dots, x_N]$.

\ms\noi Let $\fn$ be the maximal ideal $\langle x_1,\dots,
x_N\rangle$. Then one has an isomorphism
$$
\begin{array}{*{3}{rcl}}
L^N & \overset\sim{\longrightarrow} & \fn/\fn^2 \\
\noalign{\ms} (b_1,\dots, b_N) & \longmapsto &
\displaystyle{\sum^N_{i=1}}b_ix_i\ \mod \fn^2 \quad.
\end{array}
$$
Now let $\fm \subset \mathcal {O}_{X,x}$ be the maximal ideal.
Then we get a surjection
$$
\varphi:\fn/\fn^2\twoheadrightarrow \fm/\fm^2\qquad ,
$$
and for a point $b=(b_1,\dots, b_N)\in L^N$ and the associated
hyperplane $H_b$, the local ring of $x$ in $X\cdot H_b$ is
$$\mathcal {O}_{X,x}/\langle\sum\limits^N_{i=1}b_i x_i\rangle.$$
By the above lemma, $x$ is an ordinary quadratic singularity if
the image $\overline{\sum b_ix_i}$ of $\sum b_ix_i$ in $\fm/\fm^2$
is non-zero, and if the associated hyperplane in
$\IP_L(\fm/\fm^2)$ intersects the quadric in $\IP_L(\fm/\fm^2)$
associated to the singularity transversally. The latter condition
defines an open subset $U'$ in the dual projective space
$\IP_L((\fm/\fm^2)^\vee)$ parametrizing the hyperplanes in
$\IP_L(\fm/\fm^2)$. Consider the non-empty open subset
$U''\subseteq\IP_L((\fn/\fn^2)^\vee)$ on which the projection
$$
p:\IP_L((\fn/\fn^2)^\vee) \dashrightarrow\IP_L((\fm/\fm^2)^\vee)
$$
associated to $\varphi^\vee$ is defined. (To wit: $U''$ is the
complement of
$\IP_L((\ker\varphi)^\vee)\subseteq\IP_L((\fn/\fn^2)^\vee)$.
Letting $U=p^{-1}(U')$, we see that for the hyperplanes $H$ in $U$
the intersection $X\cdot H$ has an ordinary quadratic singularity
at $x$.

\ms\noi {\bf (2.4)} We can now finish the proof of Theorem 2.
Applying Theorem 4 to $X_F$ and combining it with the result on
$X_K$, we obtain the wanted Lefschetz pencil over $\Spec(A)$
provided there is an $A$-rational point in $\mathcal G$,
corresponding to a line $L$ over $A$, such that $L_\eta$ lies in
the open $W_{X_\eta} \subset \mathcal G_\eta$ (constructed in
(2.1)) and $L_s$ lies in the open $U_{X_s} \subset \mathcal G_s$
(constructed in Theorem 4). This existence, under the conditions
of Theorem 2, follows now by applying the arguments in the proof
of Theorem 0 to $P=\mathcal G$, $V_1 = W_{X_\eta}$ and $V_2 =
U_{X_s}$. Note that the specialization map $\mathcal G(K) =
\mathcal G(A) \rightarrow \mathcal G(F)$ is surjective, and that
$\mathcal G_L$, over a field $L$, has a cellular decomposition, so
that $\mathcal G(L)$ is dense in $\mathcal G$ for infinite $L$.

\bigskip
\setcounter{section}{3}
\begin{center}{\bf\arabic{section}. Desingularization of ordinary quadratic singularities}\end{center}

For the applications, it is important to have a good description
of varieties with ordinary quadratic reduction, and also a description of
their desingularization, because such schemes may be non-regular.
We recall the following description of local rings around an
ordinary quadratic singularity [SGA 7, XV, 1.32].

\ms\noi {\bf Lemma 4} ~{\it Let $X$ be a flat scheme of finite
type over $A$, and assume that $X$ is smooth over $A$ except for
one singular point $x \in X_s$ which is an ordinary quadratic
singularity (in $X_s$). Assume that $X_s$ is of dimension $n$ at
$x$. Then, after possibly passing to a finite \'etale extension of
$A$, the Henselization of $\mathcal O_{X,x}$ is isomorphic to the
Henselization of the following ring $B$ at the maximal ideal
$\langle x_1,\ldots, x_{n+1}, \pi \rangle$.

(i) If $x$ is non-degenerate:
$$
B = A[x_1,\ldots, x_{n+1}]/\langle Q(x_1,\ldots,x_{n+1}) -
c\rangle,
$$
where $Q$ is a non-degenerate quadratic form over $A$ and $c \in
\fm \setminus \{0\}$.

(ii) If $x$ is degenerate (which can only happen if $char(F) = 2$
and $n=2m$ is even):
$$
B = A[x_1,\ldots,x_{n+1}]/\langle P(x_1,\ldots, x_{2m}) +
x_{n+1}^2 + bx_{n+1} + c\rangle,
$$
where $P$ is a non-degenerate quadratic form over $A$ and $b,c \in
A$ with $b^2-4c \in \fm\setminus \{0\}$.}

\ms In the situation of Lemma 4 (i), let $r = v(c)$, where $v$ is
the normalized valuation of $K$, so that $c = \eta \pi^r$, where
$\pi$ is a prime element in $A$ and $\eta$ is a unit in $A$. Then
we say that $X$ has an ordinary quadratic singularity of order
$r$. By possibly passing to a ramified extension of degree 2
(extracting a square root of $\pi$), we may assume that $r$ is
even.

In the situation of Lemma 4 (ii), by possibly passing to a
ramified extension of degree 2 (the splitting field of $x^2 + bx +
c$), and by a coordinate transformation, we may assume that $c=0$.
In this case we let $q=v(b)$, so that $b=\epsilon\pi^q$ with a
unit $\epsilon$, and say that $X$ has an ordinary quadratic
singularity of order $q$.

\ms\noi {\bf Theorem 5} ~{\it Let $X$ be as in Lemma 4, and let
$\varphi : \tilde X \rightarrow X$ be the blowing up of $X$ at the
singular point $x$. Assume $r$ is even in case 2, and $c=0$ in
case (ii). Then the strict transform $\tilde Y$ of $Y = X_s$ is
smooth, and the exceptional fibre $F_x =\varphi^{-1}(x)$ contains
a point $\tilde x \notin \tilde Y$ such that the following holds:

(a) $F_x \setminus \{\tilde x\}$ is smooth, and $\tilde Y$ and
$F_x$ intersect transversally, i.e., the scheme-theoretic
intersection of these inside $\tilde X$ is smooth.

(b) $\tilde X \setminus \{\tilde x\}$ is regular and has strict
semi-stable reduction.

(c) In case (i), if $x$ is of order $r$, then the behavior of
$\tilde X$ at $\tilde x$ is as follows: If $r > 2$, then $\tilde
x$ is ordinary quadratic of order $r-2$. If $r=2$, then $\tilde X$
is also smooth at $\tilde x$, and hence has strict semi-stable
reduction.

(d) In case (ii), if $x$ is of order $q$, then the behavior of
$\tilde X$ at $\tilde x$ is as follows: If $q > 1$, then $\tilde
x$ is ordinary quadratic of order $q-1$. If $q=1$, then $\tilde X$
is also smooth at $\tilde x$, and hence has strict semi-stable
reduction.}

\ms\noi {\bf Proof} Since blowing-ups are compatible with flat
base change, and since smoothness and type of the quadratic
singularity just depend on the Henselization of the local ring, we
may consider the rings $B$ in Lemma 4.

\ms\noi {\bf Case (i):} 1) Here the blowing-up of $B$ at the ideal
$\fn = \langle x_1,\ldots, x_{n+1},\pi\rangle$ is $\Proj(C)$, for
the $B$-algebra
$$
C = B[U_1,\ldots, U_{n+1},T]/I
$$
$$
I = \langle ~x_iU_j-x_jU_i, ~~x_iT-\pi U_i, ~~Q(U_1,\ldots,
U_{n+1})-\eta\pi^{r-2}T^2 ~\rangle,
$$
which is graded as quotient of the polynomial ring over $B$. In
fact, the coordinate ring of the affine chart $\{U_{n+1}\neq 0\}$
is
$$
A[u_1,\ldots,u_n,x_{n+1},t]/\langle
~Q(u_1,\ldots,u_n,1)-\eta\pi^{r-2}t^2,~~x_{n+1}t-\pi~\rangle,
$$
with $x_{n+1}u_i=x_i\quad(i=1,\ldots,n)$. A similar description
holds for the other charts $\{U_i\neq 0\}$. The coordinate ring
for the chart $\{T\neq 0\}$ is
$$
A[u_1,\ldots,u_{n+1}]/\langle
~Q(u_1,\ldots,u_{n+1})-\eta\pi^{r-2}~\rangle,
$$
with $\pi u_i=x_i\quad(i=1,\ldots,n+1)$. This shows that the
inverse image of $\fn$ is an invertible ideal: it is generated by
one element (by $x_{n+1}$, $x_i$, and $\pi$, respectively), which
is not a zero divisor. Moreover, the morphism $\Proj(C)
\rightarrow \Spec(B)$ becomes an isomorphism after inverting any
of the elements $x_1,\ldots , x_{n+1}, \pi$. Finally there is a
surjection of graded $B$-algebras
$$
C \longrightarrow \mathop{\oplus}\limits_{n\geq 0} ~~\fn^n,
$$
by sending $U_i$ and $T$ to $x_i$ and $\pi$ in $\fn$,
respectively. Thus, by lemma 5 below, $\Proj(C)$ is isomorphic to
$\Proj(\oplus_{n\geq 0}~ \fn^n)$, the blowing-up of $B$ in $\fn$.

\ms\noi 2) Assume $r>2$. We consider the special fiber of the
blowing-up, obtained by setting $\pi=0$. Thus its chart
$\{U_{n+1}\neq 0\}$ is
$$
\begin{array}{rl}
&\Spec(~k[u_1,\ldots,u_n,x_{n+1},t]/\langle
~Q(u_1,\ldots,u_n,1),~~x_{n+1}t~\rangle~)\\
= & \Spec(~R[x_{n+1},t]/\langle ~x_{n+1}t~\rangle~)
\end{array}
$$
where
$R=k[u_1,\ldots,u_n]/\langle~Q(u_1,\ldots,u_n,1)~\rangle.$\quad It
is reduced, with two smooth irreducible components intersecting
transversally - the first one being the locus $\{t=0\}$, the
second one being the locus $\{x_{n+1}=0\}$. A similar result holds
for the other charts $\{U_i\neq 0\}$. The chart $\{T \neq 0\}$ is
$$
\Spec(~k[u_1,\ldots,u_{n+1}]/\langle~Q(u_1,\ldots,u_{n+1})~\rangle~),
$$
which is smooth except for one ordinary quadratic singularity at
$u=(0,\ldots,0)$.

We may identify the irreducible components as follows. The strict
transform $\tilde Y$ of the special fiber of $Spec(B)$ is obtained
by blowing up
$$
\bar B = B/\langle \pi\rangle = F[x_1,\ldots,x_{n+1}]/\langle
Q(x_1,\ldots,x_{n+1})\rangle
$$
in the ideal $\bar\fn=\langle x_1,\ldots,x_{n+1}\rangle$. This is
$\Proj(\bar C)$, for
$$
\bar C = \bar B[U_1,\ldots,U_{n+1}]/\langle ~x_iU_j-x_jU_i,~
Q(U_1,\ldots,U_{n+1})~\rangle.
$$
The affine ring of the chart $\{U_i\neq 0\}$ is
$$
F[x_i,u_1,\ldots,\check{u_i},\ldots,u_{n+1}]/\langle
~Q(u_1,\ldots,1,\ldots,u_{n+1})~\rangle,
$$
where $x_iu_j=x_j\quad(j\neq i)$, $\check{u_i}$ means omission of
$u_i$, and the 1 in $Q(u_1,\ldots,1,\ldots,u_{n+1})$ is at
position $i$. This is smooth over $F$, and corresponds to the
locus $T=0$ in $\tilde X$.

The exceptional fibre $F_x$ is obtained by letting
$x_1=\ldots=x_{n+1}=0=\pi$ in $C$. For $r>2$ we get
$$
\Proj(F[U_1\ldots,U_{n+1},T]/\langle
~Q(U_1,\ldots,U_{n+1})~\rangle).
$$
In the chart $\{U_i\neq 0\}$ this corresponds to the locus
$x_i=0=\pi$ which is
$$
\Spec(F[u_1,\ldots,\check{u_i},\ldots,u_{n+1},t]/\langle
~Q(u_1,\ldots,1,\ldots,u_{n+1})~\rangle)
$$
and thus smooth. In the chart $\{T\neq 0\}$ we get
$$
\Spec(F[u_1,\ldots,u_{n+1}]/\langle Q(u_1,\ldots,u_{n+1}\rangle).
$$
This shows that the exceptional fiber has one ordinary quadratic
singular point which does not lie on $\tilde Y$. From the previous
description of the chart $\{T\neq 0\}$ for the whole blowing-up we
see that the order of the quadratic singularity is $r-2$.

\ms\noi 3) Now let $r=2$. Then the chart $\{U_{n+1}\neq 0\}$ of
the whole blowing-up is
$$
\Spec(~S/\langle ~x_{n+1}~t - \pi~\rangle~),
$$
where $S=A[u_1,\ldots,u_n,x_{n+1},t]/\langle ~Q(u_1,\ldots,u_n,1)
- \eta t^2 ~\rangle~$ is smooth over $A$ and $x_{n+1}, t$ are part
of a local parameter system where they vanish. Thus we get a
regular scheme with semi-stable reduction over $A$. The same holds
in the other charts $\{U_i\neq 0\}$. In the chart $\{T\neq 0\}$ we
get the smooth $A$-scheme
$$
\Spec(~A[u_1,\ldots,u_{n+1}]/\langle~Q(u_1,\ldots,u_{n+1}) -
\eta~\rangle~).
$$
The strict transform $\tilde Y$ of $Y$ has exactly the same
description as before; it is smooth, and it is again the locus
where $T=0$. The exceptional fiber is
$$
\Proj(~F[U_1,\ldots,U_{n+1},T]/\langle~Q(U_1,\ldots,U_{n+1}) -
\eta T^2~\rangle~)
$$
which is smooth as well. Therefore $\tilde X$ has strict
semistable reduction.

\ms\noi {\bf Case (ii):} Here the blowing-up of
$$
B= A[x_1,\ldots,x_{n+1}]/\langle ~P(x_1,\ldots, x_{2m}) +
x_{n+1}^2 + bx_{n+1}~\rangle
$$
($b\in\fm\setminus \{0\}$) in the ideal $\fn = \langle
~x_1,\ldots,x_{n+1},\pi ~\rangle$ is $\Proj(C)$, for
$$
C = B[U_1,\ldots, U_{n+1},T]/I,
$$
where the ideal $I$ is generated by the elements
$$
\begin{array}{l}
x_iU_j-x_jU_i \quad\mbox{for}\quad i,j \in \{1,\ldots,n+1\} \\
x_iT = \pi U_i \quad\mbox{for}\quad i \in \{1,\ldots,n+1\} \\
P(U_1,\ldots, U_n) + U_{n+1}^2 + \epsilon \pi^{q-1}TU_{n+1}.
\end{array}
$$
The coordinate ring of the chart $\{U_{n+1}\neq 0\}$ is
$$
A[u_1,\ldots,u_n,x_{n+1},t]/\langle
~x_{n+1}t-\pi,~~P(u_1,\ldots,u_n) + 1 + \epsilon \pi^{q-1}t
~~\rangle~).
$$
For $i \in \{1,\ldots,n\}$, the chart $\{U_i\neq 0\}$ is
$$
\Spec(~A[u_1,\ldots,u_{i-1},x_i,u_{i+1},\ldots,u_{n+1},t]/J~),
$$
where
$$
J = \langle ~x_it -
\pi,~~P(u_1,\ldots,u_{i-1},1,u_{i+1},\ldots,u_n) + u_{n+1}^2 +
\epsilon\pi^{q-1}tu_{n+1} ~\rangle,
$$
The affine ring for the chart $\{T\neq 0\}$ is
$$
A[u_1,\ldots,u_{n+1}]/\langle ~P(u_1,\ldots,u_n) + u_{n+1}^2 +
\epsilon\pi^{q-1}u_{n+1}~\rangle.
$$

The strict transform of the special fiber $X_s$ is the locus
$T=0$, and it has exactly the same description as in case (i)2),
except that the quadratic form is now $Q(U_1,\ldots,U_{n+1}) =
P(U_1,\ldots,U_n) + U_{n+1}^2$. Thus it is smooth.

The exceptional fibre corresponds to the locus $\{x_i = 0 =
\pi\}$. In the chart $\{U_i\neq 0\}$, for $i\neq n+1$, we get the
subscheme
$$
\Spec(~F[u_1,\ldots,\check{u_i},\ldots,u_{n+1},t]/\langle ~f_i~
\rangle~)
$$
with
$$
f_i = P(u_1,\ldots,1,\ldots,u_n)+u_{n+1}^2 +
\overline{\epsilon\pi^{q-1}}tu_{n+1}
$$
where $\overline{a}=a \mbox{ mod } \pi$. In the chart
$\{U_{n+1}\neq 0\}$ we get the subscheme
$$
\Spec(~F[u_1,\ldots,u_n,t]/\langle~f_{n+1}~\rangle~),
$$
$$
f_{n+1} = P(u_1,\ldots,u_n) + 1 + \overline{\epsilon\pi^{q-1}}t.
$$
These are smooth. In the chart $\{T\neq 0\}$ we get the scheme
$$
\Spec(~F[u_1,\ldots,u_{n+1}]/\langle~g~\rangle~),
$$
$$
g = P(u_1,\ldots,u_n) + u_{n+1}^2 +
\overline{\epsilon\pi^{q-1}}u_{n+1}.
$$
If $q>1$, this has one quadratic singularity of order $q-1$. If
$q=1$, the scheme is smooth, since $\partial g/\partial u_{n+1} =
\overline{\epsilon}\neq 0$. It is also clear that the strict
transform of the special fibre and the exceptional fibre intersect
transversally (in their smooth loci). Hence the claim follows.

\ms\noi {\bf Lemma 5} ~{\it Let $B$ be a noetherian ring, let
$I\subset B$ be an ideal, and let $\tilde X = \Proj(\oplus~I^n)$
be the blowing-up of $X = \Spec(B)$ in the closed subscheme $Y =
\Spec(B/I)$ corresponding to $I$. Let
$$
\varphi:~C \longrightarrow \bigoplus_{n\geq 0}~I^n
$$
be a surjection of graded $B$-algebras. Then the $X$-morphism
$$
f = \varphi^\ast:~ \Proj(\oplus~I^n) = \tilde X~ \longrightarrow
~Z = \Proj(C)
$$
induced by $\varphi$ is an isomorphism if and only if the
following two conditions hold.
\begin{itemize}
\item[(i)] $I$ generates an invertible ideal in $Z = \Proj(C)$.
\item[(ii)] $g$ induces an isomorphism $g^{-1}(X\setminus Y)
\mathop{\rightarrow}\limits^{\sim} X\setminus Y$.
\end{itemize}
}

\ms\noi {\bf Proof:} The two conditions are known to hold for $Z =
\tilde X$, and by the surjectivity of $\varphi$, the morphism $f$
is a closed immersion. In particular, $f$ is affine. By (i), each
point in $Z$ has an open affine neighborhood $V = \Spec(R)
\subset Z$ over which the image of $I$ is generated by one element
$a \in R$ which is not a zero divisor. Hence $f^{-1}(V) \to V$
corresponds to a surjection of rings $R \to R/J$, which induces an
isomorphism after inverting $a$, by condition (ii) (for $Z$ and
$\tilde X$). This means that $J_a = 0$ for the localization of the
ideal $J$ with respect to $a$. It follows that $J=0$, because $a$
is not a zero divisor. Therefore $f$ is an isomorphism.

\bigskip
\setcounter{section}{4}
\begin{center}{\bf\arabic{section}. Application to class field theory
of varieties over local fields with (almost) good reduction}\end{center}

Now assume that $A$ is a Henselian discrete valuation ring with
finite residue field $F = A/\mathfrak m$. Thus the fraction field
$K$ of $A$ is a non-archimedean local field (in the usual sense if
$A$ is complete). Let $V$ be a proper variety over $K$. Then we
have the \it reciprocity map for $V$ \rm
$$
\recmap V \scs  SK_1(V) \to \piab V
$$
introduced in [Bl], [Sa1] and [KS1]. Here $\piab V$ is the
abelianized algebraic fundamental group of $V$ and
$$
SK_1(V)=\Coker(\sumd V 1 K_2(k(y)) \rmapo{\partial} \sumd V 0
K_1(k(x)))
$$
where $V_i$ denotes the set of points $x \in V$ of dimension $i$,
$K_q(k(x))$ denotes the $q$-th algebraic $K$-group of the residue
field $k(x)$ of $x$, and $\partial$ is induced by tame symbols.
For an integer $n>0$ prime to $\ch(K)$ let
$$
\recmapn V \scs  SK_1(V)/n \to \piab V/n
$$
denote the induced map. Finally, for a field $L$ and a prime
$\ell$ invertible in $L$ recall the following

\medskip\noi
{\bf Conjecture $\mathbf BK_q(L,\ell)$} The Galois cohomology
group $H^q(L,\mathbb Q_\ell/\mathbb Z_\ell(q))$ is divisible.

\medskip\noi
Here and below $(q)$ denotes the usual $q$-fold Tate twist. This
conjecture is a consequence of a conjecture of Bloch and Kato
asserting the surjectivity of the symbol map $K^M_q(L) \rightarrow
H^q(L,\mathbb Z/\ell\mathbb Z(q))$ from Milnor K-theory to Galois
cohomology. The above form is weaker if restricted to particular
fields $L$, but known to be equivalent if stated for all fields.
By Kummer theory, $\BKl 1 L$ holds for any $L$ and any $\ell$. The
celebrated work of Merkuriev and Suslin [MS] shows that $\BKl 2 L$
holds for any $L$ and any $\ell$. Voevodsky [V1] proved $\BKtwo q L$
for any $L$ and any $q$. It has been announced by Rost \cite{Ro} and
Voevodsky \cite{V2} (see also \cite{SJ}) that it it holds in general.

\medskip\noi

{\bf Theorem 6} \it Let $V$ be a connected smooth projective variety
over $K$ of dimension $\geq 2$ with good or ordinary quadratic
reduction, and let $\ell$ be a prime invertible in $K$.
Assume $\BKl 3 {k(Z)}$ for any proper smooth surface $Z$
lying on $V$ which has almost good reduction over
$\Spec(\cO_{K'})$ for some finite extension $K'$ of $K$. Then
$\recmapln V$ is an isomorphism for all $\nu>0$.

\rm
\proof Let $n=\ell^\nu$, and let $X$ be as in Definition 2. By
[JS (1-6) and (7-1)] there exists a fundamental exact sequence
$$
H^K_2(V,\mathbb Z/n\mathbb Z) \to SK_1(V)/n \rmapo{\recmapn V} \piab V/n
\to H^K_1(V,\mathbb Z/n\mathbb Z) \to 0 \leqno(\mbox{4.1})
$$
where $H^K_i(V,\mathbb Z/n\mathbb Z)$ denotes the $i$-th Kato
homology of $V$ with $\mathbb Z/n\mathbb Z$-coefficients
introduced in [JS Def. 1.2]. In order to show the
surjectivity of $\recmapn V$, it thus suffices to show
$H^K_1(V,\mathbb Z/n\mathbb Z)=0$.
According to Theorem 5,
by blowing up $X$ at the points where $X_s$ is not smooth, and repeating
this finitely many times, we obtain a scheme $\tX$ which is regular with
strict semi-stable reduction, and has the same generic fibre $V$, by construction.

Then, by the isomorphism
$H^K_1(V,\mathbb Z/n\mathbb Z) \cong H^K_1(\tX_s,\mathbb Z/n\mathbb Z)$
proved in [JS Thm. 1.5 (2)] it suffices to show the
vanishing of the latter group, which is the Kato homology of the
special fibre $\tX_s$ as defined in [JS Def. 1.2].
Moreover it is easy to see that the configuration
complex $\Gamma_{\tX_s}$ of $\tX_s$ introduced in [JS Rem. 3.8]
is contractible. Thus the isomorphism
$H^K_1(\tX_s,\mathbb Z/n\mathbb Z) \cong H_1(\Gamma_{\tX_s},\mathbb Z/n\mathbb Z)$
proved in [JS Thm. 1.4] gives the desired vanishing.

Now we show the injectivity of $\recmapn V$. Let $X
\subset \Bbb P^N_A$ be an embedding of $X$ into the projective
space over $A$, and fix a prime $q$. By Theorem 2, after possibly
taking the base change with a finite unramified covering of $A$ of
$q$-power degree, there exists a Lefschetz pencil $\{V_t\}_{t\in
D}$, where $D \cong \Bbb P^1$ is a $K$-line in the dual projective
space of $\Bbb P^N_K$, satisfying the following conditions:
\begin{description}
\item[(1)] The axis of the pencil has good reduction over $K$.
\item[(2)] There exists a finite subset $\Sigma\subset \bP^1_K$ of
closed points such that for every $t\not\in \Sigma$, $V_t$ has
ordinary quadratic reduction over $\Spec(\cO_{k(t)})$.
\end{description}
\rm\medbreak

We write $\VS=\cup_{t\in \Sigma} V_t$. Fix an integer $n>0$.

\Claim {1} \it Any element $\alpha\in SK_1(V)/n$ is represented by
\begin{description}
\item[$(*)$] $ \sum_{i=1}^r f_i\otimes [x_i] \quad\text{ with }
f_i\in k(x_i)^*,$
\end{description}
where $x_i$ is a closed point of $V\backslash \VS$ for $1\leq
i\leq r$. \rm

\proof By definition $\alpha$ is represented by a sum of the form
$(*)$, where however the $x_i$ may lie on $\VS$. By a standard
Bertini argument there exists a proper smooth curve $C$ over $k$
lying on $V$ such that
\begin{description}
\item[(1)] $C$ is not contained in any fiber $V_t$. \item[(2)]
$\{x_1,\ldots,x_n\} \subset C$.
\end{description}

By (1) $\CS:=C\cap \VS$ is finite and we put $U=C\backslash \CS$.
By (2) $\alpha$ lies in the image of $SK_1(C)\to SK_1(V)$. The
claim follows from the surjectivity of the natural map
$$ \bigoplus_{x\in U_0} k(x)^* \to SK_1(C)/n \leqno(\mbox{7-4})$$
which is a consequence of the class field theory of curves over
local fields. Indeed, by [Sa1] one knows that the reciprocity map
$SK_1(C)/n\to \piab C/n$ is injective and that every character of
$\piab C/n$ which is trivial on the image of (7-4) is trivial on
$SK_1(C)/n$. Namely, if every closed point of $U$ splits completely
in a given abelian covering of $C$, each point of $\CS$ splits
completely as well. \medbreak

Now, fixing $n=\ell^\nu$, we show the injectivity of $\recmapn V$
by induction on $d:=\dim(V)$. The case $d=2$ follows from [JS].
In fact, by (4.1) it suffices to show $H^K_2(V,\mathbb Z/n\mathbb Z) =0$.
But for $d=2$ we have $H^K_3(V,\mathbb Q_\ell/\mathbb Z_\ell) = 0$
by definition, and [JS Lem.7.3] gives an isomorphism
$H^K_2(V,\mathbb Z/n) \cong H^K_2(V,\mathbb Q_\ell/\mathbb Z_\ell)[n]$,
where $B[n]=\{b\in B\mid nb=0\}$ for an abelian group $B$.
Moreover, by [JS Thm.1.6] we have an isomorphism
$H^K_2(V,\mathbb Q_\ell/\mathbb Z_\ell) \cong H^K_2(\tX_s,\mathbb Q_\ell/\mathbb Z_\ell)$,
and by [JS Thm.1.4] we have an isomorphism $H^K_2(\tX_s,\mathbb Q_\ell/\mathbb Z_\ell)
\cong H_2(\Gamma_{\tX_s},\mathbb Q_\ell/\mathbb Z_\ell)$. Since the latter
group vanishes by contractibility of $\Gamma_{\tX_s}$, the claim follows.

Now assume $d>2$. Let $\alpha\in SK_1(V)/n$ and
assume $\recmapn V(\alpha)=0$. We want to show $\alpha=0$.
Take a Lefschetz pencil as in Theorem 2. By
Claim (1) there exist $t_1,\dots,t_m\in \bP^1_k$ such that
$Y_i:=V_{t_i}$ has ordinary quadratic reduction over $k$ and
$\alpha$ lies in the image of
$$ SK_1(Y)/n \to SK_1(V)/n \quad \hbox{ with }Y:=\cup_{i=1}^{\,m} Y_i.$$
We have the commutative diagram
$$ \begin{matrix}
SK_1(Y)/n &\rmapo{\recmapn Y}& \piab Y /n \\
\downarrow&&\downarrow\\
SK_1(V)/n &\rmapo{\recmapn V}& \piab V /n \\
\end{matrix}$$

\Claim {2} \it The right vertical map is injective.

\proof \rm By Poincar\'e duality, we have isomorphisms
$$\piab V/n\isom H^{2d+1}(V,\nz(d+1)),\quad
\piab Y/n\isom H^{2d+1}_Y(V,\nz(d+1)).$$ But in the localization
sequence
$$ H^{2d}(V\backslash Y,\nz(d+1))\to H^{2d+1}_Y(V,\nz(d+1))\to
H^{2d+1}(V,\nz(d+1))$$ we have $H^{2d}(V\backslash Y,\nz(d+1))=0$
since $V\backslash Y$ is affine and $2d>\dim(V)+2=d+2$ by the
assumption that $d>2$. This proves the claim. \medbreak

By Claim (2) the desired assertion follows if we show that
$\recmapn Y$ is an isomorphism. Let $A$ be the axis of the pencil.
We have a commutative diagram
$$ \begin{matrix}
\bigoplus_{i=1}^{m-1} SK_1(A)/n &\isom& \bigoplus_{1=1}^{m-1}\piab {A} /n \\
\downarrow&&\downarrow\\
\bigoplus_{i=1}^{\,m} SK_1(Y_i)/n &\isom& \bigoplus_{i=1}^{\,m}\piab {Y_i} /n \\
\downarrow&&\downarrow\\
SK_1(Y)/n &\rmapo{\recmapn Y}& \piab Y /n \\
\downarrow\\
0\\
\end{matrix}$$
Here the vertical maps from the second to the third line are 
simply induced by the inclusions $\rho_i: Y_i \hookrightarrow Y$.
The vertical maps above are obtained by noting that $A = Y_1\cap Y_2 = 
Y_2\cap Y_3 = \ldots = Y_{m-1}\cap Y_m$, and using the inclusions
$Y_i\cap Y_{i+1} \hookrightarrow Y_i$ and $Y_i\cap Y_{i+1} \hookrightarrow Y_{i+1}$.
Then the exactness of the left vertical sequence is easily seen.
Similarly, one can show that the right vertical sequence is exact,
by using the Poincar\'e duality and the localization sequence of
\'etale cohomology. By induction hypothesis, the two upper
horizontal maps are isomorphisms. This shows that $\recmapn Y$ is
an isomorphism and completes the proof of Theorem 6.

\bigskip

    Uwe Jannsen,
    Fakult\"{a}t f\"ur Mathematik,
    Universit\"at Regensburg,
    Universit\"atsstr. 31,
    93040 Regensburg,
    GERMANY,
    uwe.jannsen@mathematik.uni-regensburg.de

\bigskip
    Shuji Saito,
    Graduate School of Mathematics,
    University of Tokyo,
    Komaba, Meguro-ku,
    Tokyo 153-8914,
    JAPAN,
    sshuji@ms.u-tokyo.ac.jp

\end{document}